\documentclass[10pt,a4paper]{article}
\usepackage[latin1]{inputenc}
\usepackage[T1]{fontenc}
\usepackage{amsmath,amscd}
\usepackage{graphicx}
\usepackage[all]{xy}
\usepackage{amsthm}
\usepackage{amssymb}
\usepackage[charter]{mathdesign}
\usepackage{amsfonts}
\usepackage{mathrsfs}
\usepackage{amsxtra}
\usepackage{bbm}
\usepackage{ae}
\usepackage[all]{xy}

\newtheorem{sats}{Theorem}[section]

\newtheorem{lem}[sats]{Lemma}

\newtheorem{prop}[sats]{Proposition}

\newcommand{\R}{\mathbbm{R}}
\newcommand{\C}{\mathbbm{C}}

\newcommand{\Z}{\mathbbm{Z}}

\newcommand{\ellL}{\mathcal{L}}

\newcommand{\id}{ \mathrm{id}}

\newcommand{\Ko}{\mathcal{K}}
\newcommand{\Bo}{\mathcal{B}}

\newcommand{\Eg}{\mathcal{E}}

\newcommand{\Hom}{\mathrm{Hom}}

\renewcommand{\epsilon}{\varepsilon}
\renewcommand{\phi}{\varphi}

\newcommand{\e}{\mathrm{e}}

\newcommand{\sla}{\!\!\!\!/\,}

\title{The Pimsner-Voiculescu sequence for coactions of compact Lie groups}
\author{Magnus Goffeng}
\date{Department of Mathematical Sciences, Division of Mathematics\\
Chalmers university of Technology and University of Gothenburg}

\begin{document}
\maketitle

\begin{abstract}
The Pimsner-Voiculescu sequence is generalized to a Pimsner-Voiculescu tower describing the $KK$-category equivariant with respect to coactions of a compact Lie group satisfying the Hodgkin condition. A dual Pimsner-Voiculescu tower is used to show that coactions of a compact Hodgkin-Lie group satisfy the Baum-Connes property. 
\end{abstract}
\normalsize

\section*{Introduction}

When $G$ is a second countable, locally compact group and $A$ is a separable $C^*$-algebra with a continuous $G$-action, the Baum-Connes conjecture states that the $K$-theory of the reduced crossed product $A\rtimes_r G$ can be calculated by means of geometric and representation theoretical properties of $G$ and $A$, see more in \cite{bch}. To be more precise, the Baum-Connes conjecture states that the assembly mapping $\mu_A:K_*^G(\Eg G;A)\to K_*(A\rtimes_r G)$ is an isomorphism. The space $\Eg G$ is the universal proper $G$-space and $K_*^G(\Eg G;A)$ is the proper equivariant $K$-homology with coefficients in $A$. There are known counterexamples when $\mu_A$ is not an isomorphism, so it is more natural to speak of groups having the Baum-Connes property. In \cite{neme}, the equivariant $K$-homology with coefficients in $A$ was proved to be the left derived functor of $F(A)=K_*(A\rtimes_r G)$ and the assembly mapping being the natural transformation from $\mathbbm{L}F$ to $F$. The approach to the Baum-Connes property using triangulated categories can be generalized to discrete quantum groups, see \cite{mehom}, which indicates that geometric techniques such as universal proper $G$-spaces can be generalized to discrete quantum groups. 

The generalization of the Baum-Connes property to quantum groups has been studied in for instance \cite{nemetva} and \cite{voigtfree}. The case studied in \cite{nemetva} is that of quantum group actions of the dual of a compact Lie group which correspond to coactions of the Lie group. In \cite{nemetva} duals of compact Lie groups were shown to satisfy the strong Baum-Connes property, i.e. the embedding of the triangulated category generated by proper coactions, the $C^*$-algebras that are Baaj-Skandalis dual to trivial $G$-actions, into the $KK$-category equivariant with respect to coactions is essentially surjective. In this paper we construct an analogue of the Pimsner-Voiculescu sequence for coactions of a compact Hodgkin-Lie group $G$ that describes how the $KK$-category equivariant with respect to coactions of $G$ is built up from the $C^*$-algebras with coactions of $G$ which are proper in the sense of \cite{nemetva}.

The starting point is to express the Pimsner-Voiculescu sequence for $\Z$-actions in terms of a property of the representation ring of a rank one torus. Using the Universal Coefficient Theorem, the Pimsner-Voiculescu sequence can be constructed from a Koszul complex 
\[0\to R(T)\xrightarrow{\alpha} R(T)\to 0,\] 
where $\alpha$ is defined as multiplication by $1-t$ under the isomorphism $R(T)\cong \Z[t,t^{-1}]$. When $A$ has a coaction of $T$, i.e. a $\Z$-action, the tensor product over $R(T)$ between this Koszul complex and $K_*^T(A\rtimes_r \Z)$ gives the Pimsner-Voiculescu sequence. In the generalization to higher rank, when $T$ is a torus of rank $n$ we consider the Koszul complex
\[0\to  \wedge^{n}R(T)^n \to\wedge^{n-1}R(T)^n\to\ldots \to\wedge^{2}R(T)^n\to R(T)^n\to R(T)\to 0.\]
The boundary mappings in this complex are defined from interior multiplication with the element $\sum (1-t_i)e_i^*\in \Hom_{R(T)}(R(T)^n,R(T))$. If $G$ is a compact Hodgkin-Lie group with maximal torus $T$, the representation ring $R(T)$ is a free $R(G)$-module by \cite{stei}, so the generalization from a torus to compact Hodgkin-Lie groups goes in a straightforward fashion. Just as when the rank is $1$, the Koszul complex above can be used to produce sequence of distinguished triangles which is the analogue of a Pimsner-Voiculescu sequence for the $K$-theory of crossed products by coactions of $G$.

We will give a geometric description of a sequence of distinguished triangles in the $KK$-category equivariant with respect to coactions of $G$ that corresponds to the above Koszul complex under the Universal Coefficient Theorem. As for the Pimsner-Voiculescu sequence for $\Z$ we will obtain a projective resolution of the crossed product by a coaction in the sense of triangulated categories rather than exact sequences. Using suitable tensor products we produce in Theorem \ref{gpv} a sequence of distinguished triangles in the $KK$-category equivariant with respect to coactions of $G$ that we call the generalized Pimsner-Voiculescu tower for $A$:
\begin{align*}
\xymatrix@C=0.2em@R=2.71em{
& \C^w\otimes A\ar[rr]&&
\Sigma^nD_{n-1}(A) \ar[rr] \ar[dl] &&
\Sigma^nD_{n-2}(A) \ar[rr] \ar[dl] &&
\cdots \\
&& \Sigma\C^{wn}\otimes A \ar[ul]|-\circ&&
\Sigma^{2}\C^{wk_{n-1}}\otimes A \ar[ul]|-\circ\ar[ll]|-\circ&&
\cdots\ar[ul]|-\circ\ar[ll]|-\circ
}\\
\xymatrix@C=0.2em@R=2.71em{
\cdots \ar[rr] &&
\Sigma^nD_2(A) \ar[dl]\ar[rr]&&
\Sigma^nD_1( A) \ar[dl]\ar[rr]&&
t(A\rtimes_r\hat{G})\ar[dl] \\
& \cdots &&
\Sigma^{n-1}\C^{wk_2}\otimes A \ar[ul]|-\circ \ar[ll]|-\circ &&
\Sigma^n\C^w\otimes A \ar[ul]|-\circ \ar[ll]|-\circ
}.
\end{align*}
Here $t(A\rtimes_r\hat{G})$ denotes the $C^*$-algebra $A\rtimes_r \hat{G}$ equipped with the trivial $\hat{G}$-action and the terms $D_i(A)$ can be explicitly described as a braided tensor product. Taking $K$-theory of the lower row will give a complex similar to the Koszul complex that in a sense forms a projective resolution of the $K$-theory of $A\rtimes \hat{G}$. The dual Pimsner-Voiculescu gives a more precise description of the results of \cite{nemetva} by a sequence of distinguished triangles in $KK^G$ that describes the crossed product $A\rtimes_r \hat{G}$ in terms of $G-C^*$-algebras with trivial $G$-action, thus giving a direct route to the strong Baum-Connes property of $\hat{G}$.

The paper is organized as follows; the first section consists of a review of $KK$-theory of actions and coactions. In particular we gather some known results about the braided tensor product and the Drinfeld double which plays a mayor role in constructing the dual Pimsner-Voiculescu tower. The main references of this section are \cite{baska}, \cite{baskauni}, \cite{bavae}, \cite{kuva}, \cite{neme}, \cite{voigt} and \cite{vaind}. In the second section a geometric construction of the Pimsner-Voiculescu sequence for $\Z$-actions is presented and generalized to higher rank via a Koszul complex. In the third section the restriction functor for coactions is used to generalize the Pimsner-Voiculescu sequence to coactions of compact Hodgkin-Lie groups $G$. As an example of this we calculate the $K$-theory of some compact homogeneous spaces. By similar methods, a dual Pimsner-Voivulescu tower is constructed in $KK^G$, following the ideas of \cite{neme}. At the end of the paper we discuss some possible generalizations to duals of Woronowicz deformations.

\paragraph{Acknowledgments} 
The author would like to thank Ryszard Nest for posing the question on how to explicitly construct the crossed product of a $\hat{G}-C^*$-algebra from trivial actions and for much inspiration in the writing process.

\section{Actions and coactions of compact groups}

The standard approach to equivariant $K$-theory is to introduce equivariant $KK$-theory. If $A$ and $B$ are two separable $C^*$-algebras with a continuous action of a locally compact group $G$, the equivariant $KK$-group $KK^G(A,B)$ is defined as the set of homotopy classes of $G$-equivariant $A-B$-Kasparov modules which forms an abelian group under direct sum. The $KK$-groups can be equipped with a product such that if $C$ is a third separable $C^*$-algebra with a continuous $G$-action there is an additive pairing called the Kasparov product
\[KK^G(A,B)\times KK^G(B,C)\to KK^G(A,C).\]
Following the standard construction, we let $KK^G$ denote the additive category of all separable $C^*$-algebras with a continuous $G$-action and a morphism in $KK^G$ from $A$ to $B$ is an element of $KK^G(A,B)$. The composition of two $KK^G$-morphisms is defined to be their Kasparov product. The group $KK^G(\C,A)$ coincides with the equivariant $K$-theory of $A$. In particular, if $G$ is compact $KK^G(\C,\C)=R(G)$, the representation ring of $G$. The action of $R(G)$ on equivariant $K$-theory generalizes to an $R(G)$-module structure on the bivariant groups $KK^G(A,B)$. 

The category $KK^G$ can be equipped with a triangulated structure with a mapping cone coming from the mapping cone construction of a $*$-homomorphism. The triangulated structure on $KK^G$ is universal in the sense that any homotopy invariant, stable, split-exact functor on the category of $C^*$-algebras with a continuous $G$-action defines a homological functor on $KK^G$. The construction of the triangulated structure and its universality are thoroughly explained in \cite{neme}. Let us just recall the basics of the construction of the triangulated structure on $KK^G$. The suspension $\Sigma A$ of a $G-C^*$-algebra is defined by $C_0(\R)\otimes A$. By Bott periodicity $\Sigma^2\cong \id$. A distinguished triangle in $KK^G$ is a triangle isomorphic to one of the form 
\[\xymatrix@C=0.8em@R=2.71em{
C(f) \ar[rr] & & A \ar[dl]^f\\
  &B, \ar[ul]|-\circ
}\]
where $C(f)$ is the mapping cone of the equivariant $*$-homomorphism $f:A\to B$. In particular, if $f:A\to B$ is a surjection and admits an equivariant completely positive splitting the natural mapping $\ker(f)\to C(f)$ defines an equivariant $KK$-isomorphism, so under suitable assumptions a distinguished triangle is isomorphic to a short exact sequence. 

How to construct $KK$-theory of coactions of groups is easiest seen in the simpler case when $G$ is an abelian group. If $A$ is a $C^*$-algebra equipped with an action $\alpha$ of the abelian group $G$, the crossed product $A\rtimes_r G$ carries a natural action of the Pontryagin dual $\hat{G}$. This action is called the dual action of $\hat{G}$. Since abelian groups are exact, the crossed product by an abelian group defines a triangulated functor $KK^G\to KK^{\hat{G}}$. The crossed product by the dual action is described by Takesaki-Takai duality which states that there is an equivariant isomorphism 
\[A\rtimes_rG\rtimes_r\hat{G}\cong A\otimes \Ko(L^2(G)),\]
where $A\rtimes_rG\rtimes_r\hat{G}$ is equipped with the dual action of $G$ and the $G$-action on $A\otimes \Ko(L^2(G))$ is defined as $\alpha\otimes Ad$. Takesaki-Takai duality implies that the crossed product defines a triangulated equivalence $KK^G\to KK^{\hat{G}}$.

An action $\alpha$ of a group $G$ on $A$ defines a $*$-homomorphism $\Delta_\alpha:A\to \mathcal{M}(A\otimes C_0(G))$ by letting $\Delta_\alpha(a)$ be the function $g\mapsto \alpha_g(a)$. When $G$ is abelian there is a natural isomorphism $C_0(\hat{G})\cong C^*_r(G)$ and a $\hat{G}$-action corresponds to a non-degenerate $*$-homomorphism $\Delta_A:A\to \mathcal{M}(A\otimes _{min}C^*_r(G))$ satisfying certain conditions. The first instance of a coaction of a group $G$ is on $C^*_r(G)$. Using the universal property of $C^*_r(G)$, one can construct a non-degenerate mapping $\Delta:C^*_r(G)\to \mathcal{M}(C^*_r(G)\otimes_{min}C^*_r(G))$ called the comultiplication and is induced from the diagonal homomorphism $G\to G\times G$. Clearly, the mapping $\Delta$ satisfies:
\[(\Delta\otimes \id)\Delta=(\id\otimes \Delta)\Delta,\]
so we say that $\Delta$ is coassociative. Since $\Delta_{21}=\Delta$ the comultiplication $\Delta$ is cocommutative, so if we interpret $C^*_r(G)$ as the functions on a reduced locally compact quantum group $\hat{G}$ then $\hat{G}$ can be thought of as abelian, see more in \cite{kuva}. With the abelian setting as motivation, we say that a separable $C^*$-algebra $A$ has a coaction of the locally compact second countable group $G$ if there is non-degenerate $*$-homomorphism $\Delta_A:A\to \mathcal{M}(A\otimes _{min}C^*_r(G))$ satisfying the two conditions that $\Delta_A(A)\cdot 1_A\otimes _{min}C^*_r(G)$ is a dense subspace of $A\otimes _{min}C^*_r(G)$ and that $\Delta_A$ is coassociative in the sense that 
\begin{equation}
\label{coass}
(\Delta_A\otimes \id_{C^*_r(G)})\Delta_A=(\id_A\otimes \Delta)\Delta_A.
\end{equation}
A separable $C^*$-algebra equipped with a coaction of $G$ will be called a $\hat{G}-C^*$-algebra. Sometimes we will abuse the notation and call a coaction of $G$ a $\hat{G}$-action. An example of a coaction is the dual coaction on $C^*$-algebras of the form $A=B\rtimes_rG$, for some $G-C^*$-algebra $B$. When $G$ is discrete we can decompose $B\rtimes_rG$ by means of the dense subspace $\oplus_{g\in G} B\lambda_g$ and the dual coaction is defined by $\Delta_A(b\lambda_g):=b\lambda_g\otimes \lambda_g$. In the general setting, the construction of the dual coaction goes analogously and we refer the reader to \cite{baska}. 

Much of the theory for group actions also hold for group coactions, the crossed product will as for abelian groups be a stepping stone back and forth between actions and coactions. In \cite{baska}, the $KK$-theory equivariant with respect to a bi-$C^*$-algebras and the corresponding Kasparov product was constructed. In \cite{voigt} it was proved that the $KK$-theory equivariant with respect to a locally compact quantum group has a triangulated structure defined in the same fashion as for a group. 

Let us explain the setting of \cite{baska} more explicitly in the case of coactions of a group. An $A-B$-Hilbert bimodule $\Eg$ is called $\hat{G}$-equivariant if there is a coaction $\delta_\Eg:\Eg\to \ellL_{B\otimes_{min}C^*_r(G)}( B\otimes_{min}C^*_r(G),\Eg\otimes C^*_r(G))$ satisfying a coassociativity condition similar to \eqref{coass} and $\delta_\Eg$ should commute with the $A$-action and $B$-action in the obvious ways. By Proposition $2.4$ of \cite{baska}, the coaction $\delta_\Eg$ is uniquely determined by a unitary $V_\Eg \in \ellL(\Eg\otimes _{\Delta_B}(B\otimes_{min}C^*_r(G)), \Eg\otimes C^*_r(G))$ via the equation $\delta_\Eg(x)y=V_\Eg(x\otimes_{\Delta_B}y)$ for $x\in \Eg$ and $y\in B\otimes_{min} C^*_r(G)$. A $\hat{G}$-equivariant $A-B$-Kasparov module is an $A-B$-Kasparov module $(\Eg,F)$ such that $\Eg$ is a $\hat{G}$-equivariant $A-B$-Hilbert module and the operator $F$ commutes with the unitary $V_\Eg$ up to a compact operator. The group $KK^{\hat{G}}(A,B)$ is defined as the homotopy classes of $\hat{G}$-equivariant $A-B$-Kasparov modules. The additive category $KK^{\hat{G}}$ is defined by taking the objects to be separable $\hat{G}-C^*$-algebras and the group of morphisms from $A$ to $B$ is $KK^{\hat{G}}(A,B)$. The composition in $KK^{\hat{G}}$ is Kasparov product of $\hat{G}$-equivariant Kasparov modules. 

To a closed subgroup $H$ of $G$, the restriction of a $G$-action to $H$ defines a restriction functor $Res^G_H:KK^G\to KK^H$ and its right adjoint is the induction functor $Ind_H^G:KK^H\to KK^G$. However the restriction goes in the other direction for coactions. When $H$ is a closed subgroup of $G$, there is a non-degenerate embedding $C^*(H)\subseteq \mathcal{M}(C^*(G))$ so a coaction of $H$ can be restricted to a coaction of $G$. This construction defines a triangulated functor $Res^{\hat{H}}_{\hat{G}}:KK^{\hat{H}}\to KK^{\hat{G}}$.

The crossed product $B\mapsto B\rtimes_rG$ sends a $G-C^*$-algebra to a $\hat{G}-C^*$-algebra and if $G$ is exact the crossed product induces a triangulated functor $KK^G\to KK^{\hat{G}}$. In order to construct a duality similar to Takesaki-Takai duality one introduces the crossed product by a coaction. If $A$ is a $\hat{G}-C^*$-algebra we define 
\[A\rtimes_r\hat{G}:=[\Delta_A(A)\cdot 1_A\otimes C_0(G)]\subseteq \mathcal{M}(A\otimes \Ko(L^2(G))).\]
It follows from Lemma $7.2$ of \cite{baskauni} that $A\rtimes_r\hat{G}$ forms a $C^*$-algebra. For a thorough introduction to crossed products by coactions see \cite{raeburn}. The $C^*$-algebra $A\rtimes_r\hat{G}$ carries a continuous $G$-action defined in the dense subspace $\Delta_A(A)\cdot 1_A\otimes C_0(G)$  by
\[g.(\Delta_A(a)\cdot 1_A\otimes f):=\Delta_A(a)\cdot 1_A\otimes g.f.\]
Similarly to the abelian setting, Takesaki-Takai duality holds so there are equivariant isomorphisms $B\rtimes_rG\rtimes_r\hat{G}\cong B\otimes \Ko(L^2(G))$ and $A\rtimes_r\hat{G}\rtimes_rG\cong A\otimes\Ko(L^2(G))$ which ensures that the crossed product defines an equivalence of triangulated categories known as Baaj-Skandalis duality.

The tensor product on the category of $G-C^*$-algebras is well defined. If $A$ and $B$ have actions $\alpha$ respectively $\beta$ of $G$ the tensor product $A\otimes_{min}B$ can be equipped with the action $\alpha\otimes \beta:G\to Aut(A\otimes_{min}B)$. However, for a non-abelian group $G$ the construction of a tensor product of $\hat{G}-C^*$-algebras can not be done by just taking tensor products of the $C^*$-algebras. The tensor product relevant for $\hat{G}-C^*$-algebras is the braided tensor product over $\hat{G}$ which requires one further structure. Suppose that $A$ is a $\hat{G}$-algebra with a continuous $G$-action $\alpha$. If the action $\alpha$ satisfies that 
\begin{equation}
\label{ydcond}
\Delta_A\circ \alpha_g=(\alpha_g\otimes Ad(g))\Delta_A
\end{equation}
we say that $A$ is a Yetter-Drinfeld algebra. An example of a Yetter-Drinfeld algebra is $C^*_r(G)$ with $G$-action defined by the adjoint action $G\to Aut(G)$. It is much easier to construct a Yetter-Drinfeld algebra from a $G-C^*$-algebra, if $A$ is a $G-C^*$-algebra we can in a functorial way define a coaction of $G$ on $A$ by setting $\Delta_A(a):=a\otimes 1$. When $A$ is a Yetter-Drinfeld algebra, the $C^*$-algebra $A\rtimes_r\hat{G}$ is also a Yetter-Drindeld algebra since the morphism $\Delta_A$ is covariant with respect to the $G$-action and $\Delta_A$ extends to a coaction of $G$ on $A\rtimes_r \hat{G}$, see more in \cite{voigt}. This construction is functorial and the crossed product can be seen as a functor on the category of Yetter-Drinfeld algebras.

When $A$ is a Yetter-Drinfeld algebra and $B$ is a $\hat{G}-C^*$-algebra we define the mappings 
\begin{align*}
\iota_A&:A\to \mathcal{M}(A\otimes_{min} B\otimes \Ko(L^2(G))),\quad \iota(a):=\Delta_\alpha(a)_{13}\\
\iota_B&:B\to \mathcal{M}(A\otimes_{min} B\otimes \Ko(L^2(G))),\quad \iota(b):=\Delta_B(b)_{23}.
\end{align*}
Following \cite{voigt}, the braided tensor product $A\boxtimes _{\hat{G}} B$ is defined as the closed linear span of $\iota_A(A)\cdot \iota_B(B)$. By Proposition $8.3$ of \cite{vaind}, $A\boxtimes_{\hat{G}} B$ forms a $*$-subalgebra of $\mathcal{M}(A\otimes_{min} B\otimes \Ko(L^2(G)))$ so the braided tensor product is a $C^*$-algebra. The coaction of $G$ on $A\boxtimes _{\hat{G}} B$ is defined by 
\[\Delta_A\boxtimes_{\hat{G}} \Delta_B(\iota_A(a)\cdot\iota_B(b)):=(\iota_A\otimes \id)(\Delta_A(a))\cdot (\iota_B\otimes \id)(\Delta_B(b)).\]
Observe that since $C^*_r(G)$ is cocommutative, the adjoint $\hat{G}$-action is trivial and a similar construction of a braided tensor product over $G$ between $G-C^*$-algebras with trivial $\hat{G}$-actions coincides with the usual tensor product. In general, the braided tensor product over $G$ does not need to coincide with the usual tensor product. By Lemma $3.5$ of \cite{voigt} there is a $G$-equivariant isomorphism
\begin{equation}
\label{crobra}
(A\boxtimes_{\hat{G}} B)\rtimes_r \hat{G}\cong (A\rtimes_r \hat{G})\boxtimes_{\hat{G}} B
\end{equation}
where the $G$-coaction on the right hand side is the trivial one on $B$. More generally, this identity holds for any quantum group and in particular also for braided tensor products over $G$. We will prove this statement in special case of braided tensor products over a compact group $G$ with $C(G)$ below in Lemma \ref{tenscross}.

If we interpret the structure of a Yetter-Drinfeld algebra as two actions of the quantum groups $G$ and $\hat{G}$ satisfying a certain cocycle relation, the cocycle defines a quantum group by means of a double crossed product such that Yetter-Drinfeld algebras are precisely the $C^*$-algebras with an action of this double crossed product. The right quantum group to look at is the Drinfeld double $D(G)$. Using the notations of quantum groups, the algebra of functions on $D(G)$ is $C_0(G,C^*_r(G))=C_0(G)\otimes C^*_r(G)$ with the obvious action and coaction of $G$. The action and coaction define a comultiplication 
\[\Delta_{D(G)}:C_0(D(G))\to \mathcal{M}(C_0(D(G))\otimes C_0(D(G)))\]
by $\Delta_{D(G)}:=\sigma_{23}Ad(W_{23})(\Delta_{C_0(G)}\otimes \Delta_{C^*_r(G)})$ where $W\in \Bo(L^2(G)\otimes L^2(G))$ is the multiplicative unitary of $G$ defined by $Wf(g_1,g_2)=f(g_1,g_1g_2)$. The comultiplication $\Delta_{D(G)}$ makes $D(G)$ into a quantum group by Theorem $5.3$ of \cite{bavae}. A Yetter-Drinfeld algebra $A$ with the action $\alpha$ and coaction $\Delta_A$  correspond to a $D(G)-C^*$-algebra by defining the $D(G)$-coaction 
\[\Delta_A^{D(G)} :=(\Delta_{\alpha}\otimes\id) \Delta_A:A\to \mathcal{M}(A\otimes_{min}C_0(D(G))),\] 
see more in Proposition $3.2$ of \cite{voigt}. Therefore we can consider the braided tensor product as a tensor product between $D(G)-C^*$-algebras and $\hat{G}-C^*$-algebras. The braided tensor product induces a biadditive functor 
\[\boxtimes_{\hat{G}}:KK^{D(G)}\times KK^{\hat{G}}\to KK^{\hat{G}}.\] 
Much of the theory of coactions can be done without introducing any quantum groups, but in order to construct the Pimsner-Voiculescu sequence for coactions of compact Hodgkin-Lie groups we will need the braided tensor product as a biadditive functor between $KK$-categories.

\section{The Pimsner-Voiculescu sequence from the viewpoint of representation rings}

In this section we will study the Pimsner-Voiculescu sequence for $\Z$ and generalize to a Pimsner-Voiculescu tower for $\Z^n$. We will use representation theory to calculate all the mappings explicitly. These calculations will in a surprisingly straightforward way give a natural route to a Pimsner-Voiculescu tower for coactions of compact Lie groups. 

Consider the evaluation mapping $l:C_0(\R)\to C_0(\Z)$. This mapping fits into a $\Z$-equivariant short exact sequence 
\begin{equation}
\label{firstpv}
0\to \Sigma C_0(\Z)\to C_0(\R)\xrightarrow{l} C_0(\Z)\to 0.
\end{equation}
The $\Z$-equivariant Dirac operator $D\sla$ on $\R$ defines a $\Z$-equivariant odd unbounded $K$-homology class, thus an element $[D\sla]\in KK^\Z(C_0(\R),\Sigma\C)$. While $\R$ is the universal proper $\Z$-space the element $[D\sla]$ is the Dirac element of $\Z$ and the strong Baum-Connes property of $\Z$ implies that $[D\sla]$ is a $KK^\Z$-isomorphism. The exact sequence \eqref{firstpv} induces a distinguished triangle in $KK^{\Z}$ which after using the isomorphism $C_0(\R)\cong \Sigma \C$ and rotation $4$ steps to the left becomes
\begin{equation}
\label{futt}
\xymatrix@C=0.2em@R=2.71em{
C_0(\Z) \ar[rr] & & C_0(\Z) \ar[dl]\\
  &\C. \ar[ul]|-\circ
}
\end{equation}
In a certain sense, the distinguished triangle \eqref{futt} captures the entire behavior of the Pimsner-Voiculescu sequence. If $A$ is a $\Z-C^*$-algebra we can apply Baaj-Skandalis duality to \eqref{futt} and tensor with $A\rtimes_r\Z$. If we apply Baaj-Skandalis duality again, we obtain a distinguished triangle in $KK^\Z$:
\[
\xymatrix@C=0.2em@R=2.71em{
A \ar[rr] & & A \ar[dl]\\
  &A\rtimes_r\Z, \ar[ul]|-\circ
}
\]
where $A\rtimes_r\Z$ is given the trivial $\Z$-action. Taking $K$-theory of this distinguished triangle gives back the classical Pimsner-Voiculescu sequence due to the following lemma:

\begin{prop}
\label{ranone}
When $T$ is a torus of rank $1$ and the element $\kappa\in KK^T(\C,\C)$ is defined using the isomorphisms $KK^T(\C,\C)\cong \Hom_{R(T)}(R(T),R(T))$ and $R(T)\cong \Z[t,t^{-1}]$ as 
\[\kappa f(t,t^{-1})=(1-t)f(t,t^{-1}),\] 
the $KK$-morphism $\kappa$ is Baaj-Skandalis dual to the $KK$-morphism $C_0(\Z)\to C_0(\Z)$ defined by \eqref{firstpv}.
\end{prop}

Observe that the $K$-theory of the exact sequence \eqref{firstpv} is described from the exact sequence:
\[0\to R(T)\xrightarrow{1-t} R(T)\to \Z\to 0,\]
by Proposition \ref{ranone}. The first terms in this exact sequence is the Koszul complex defined by $1-t\in \Hom_{R(T)}(R(T),R(T))$ and $\Z$ is the cohomology of the Koszul complex. 

\begin{proof}
Let $\kappa_0\in \Hom_{R(T)}(R(T),R(T))$ denote the Baaj-Skandalis dual of the $KK$-morphism induced from \eqref{firstpv}.  It follows directly from the construction that the mapping $R(T)\to \Z$ induced from $\Sigma C_0(\Z)\to C_0(\R)$ is the augmentation mapping $\Z[t,t^{-1}]\to \Z$ onto the generator of $K_1(C_0(\R))$. Therefore the image of $\kappa_0$ is the ideal generated by either $1+t$ of $1-t$ so $\kappa_0$ is of the form $u\cdot(1\pm t)$ for some unit $u\in \Z[t,t^{-1}]$. The sign and $u=1$ is found by either a direct calculation or by considering the Pimsner-Voiculescu sequence for $C_0(\Z)$.
\end{proof}

We will return to the Koszul complexes later on. First we will construct a geometric interpretation of the higher rank situation. Assume that $T$ is a torus of rank $n$ and consider the semi-open unit cube $I=[0,1[^n\subseteq \R^n$. For $i=1,\ldots, n$ we define $\tilde{X}_i$ as the set of open $i-1$-dimensional faces of $I$. The union satisfies 
\[\cup_{i=1}^n \tilde{X}_i= \partial I\cap I.\] 
We let $k_i$, for $i=1,2,\ldots n$, denote the integers  
\[k_i:=\begin{pmatrix} n\\ i-1\end{pmatrix}.\]
The set $\tilde{X}_i$ has $k_i$ connected components so if we choose a homeomorphism $]0,1[\cong \R$ there are homeomorphisms 
\begin{equation}
\label{facehom}
\tilde{X}_i\cong \coprod _{j=1}^{k_i} \R^{i-1} \quad \mbox{for}\quad i=1,2,\ldots, n,
\end{equation}
where we interpret $\R^0$ as the one-point space. We take $X_i$ to be the $\Z^n$-translates of $\cup_{j\leq i}\tilde{X}_j$ and define $Y_i:=\R^n\setminus X_i$ for $1=1,2,\ldots, n$ and $Y_0:=\R^n$. 

\begin{prop}
\label{facequo}
For $i=1,2,\ldots, n$ there are $\Z^n$-equivariant isomorphisms 
\[C_0(Y_{i-1})/C_0(Y_{i})\cong \C^{k_i}\otimes \Sigma^{i-1}C_0(\Z^n).\]
\end{prop}

\begin{proof}
By equation \eqref{facehom} there is a $\Z^n$-equivariant homeomorphism 
\[Y_{i-1}\setminus Y_i \cong \coprod_{m\in \Z^n}\left(\coprod _{j=1}^{k_i} \R^{i-1}\right),\]
where $\Z^n$ acts by translation on the first disjoint union. Therefore 
\[C_0(Y_{i-1})/C_0(Y_{i})\cong C_0(Y_{i-1}\setminus Y_{i})\cong C_0\left(\coprod_{m\in \Z^n}\left(\coprod _{j=1}^{k_i} \R^{i-1}\right)\right) \cong\]
\[\cong\C^{k_i}\otimes C_0(\Z^n\times \R^{i-1})\equiv  \C^{k_i}\otimes \Sigma^{i-1}C_0(\Z^n).\]
\end{proof}

Consider the classifying space $\R^n$ for proper actions of $\Z^n$. Since $\Z^n$ has the strong Baum-Connes property, the Dirac element $[D\sla]$ induces a $KK^{\Z^n}$-isomorphism $C_0(\R^n)\cong \Sigma^n \C$. An alternative approach to constructing this isomorphism is the Julg theorem which implies that for any $T-C^*$-algebra $A$ there is an isomorphism $K^T_*(A)\cong K_*(A\rtimes_r T)$. Therefore $K^T_*(\Sigma^n\C\rtimes \Z^n)\cong K_*^T(C_0(\R^n)\rtimes \Z^n)$ and the statement follows from the Universal Coefficient Theorem for the compact Hodgkin-Lie group $T$, see more in \cite{rosch}. 

For $i=1,2,\ldots , n$, Proposition \ref{facequo} implies that there is a $\Z^n$-equivariant short exact sequence 
\begin{equation}
\label{shortdef}
0\to C_0(Y_{i})\to C_0(Y_{i-1})\to \C^{k_i}\otimes \Sigma^{i-1}C_0(\Z^n)\to 0.
\end{equation}
We will by $\kappa_i\in KK^{\Z^n}(\C^{k_i}\otimes C_0(\Z^n),\C^{k_{i+1}}\otimes C_0(\Z^n))$ denote the $\Z^n$-equivariant $KK$-morphism defined in such a way that the extension class defined by \eqref{shortdef} composed with the restriction mapping $C_0(Y_i)\to \C^{k_{i+1}}\otimes \Sigma^iC_0(\Z^n)$ coincides with $\Sigma^{i-1}\kappa_i$. Notice that $Y_n=\Z^n\times]0,1[^n$ and $Y_0=\R^n$ so we have that $C_0(Y_{n})=\Sigma^nC_0(\Z^n)$ and $C_0(Y_0)=C_0(\R^n)$, the latter being $KK^{\Z^n}$-isomorphic to $\Sigma^n\C$. Thus we get a sequence of distinguished triangles in $KK^{\Z^n}$:
\begin{align}
\label{zntower}
\xymatrix@C=0.2em@R=2.71em{
& \Sigma^nC_0(\Z^n) \ar[rr]&&
C_0(Y_{n-1}) \ar[rr] \ar[dl] &&
\cdots \ar[dl]\\
&& \C^{n}\otimes \Sigma^{n-1}C_0(\Z^n) \ar[ul]|-\circ_{\Sigma^n\kappa_n}&&
\C^{k_{n-1}}\otimes \Sigma^{n-2} C_0(\Z^n) \ar[ul]|-\circ \ar[ll]|-\circ_{\Sigma^{n-1}\kappa_{n-1}}&&
\quad \cdots\quad \ar[ul]|-\circ \ar[ll]|-\circ}\\
\xymatrix@C=0.2em@R=2.71em{
\cdots \ar[rr] &&
C_0(Y_2) \ar[rr] \ar[dl] &&
C_0(Y_1) \ar[dl]\ar[rr]&&
\Sigma^n\C \ar[dl] \\
& \cdots \quad&& 
\C^{n}\otimes \Sigma C_0(\Z^n)\ar[ll]|-\circ_{\Sigma^2\kappa_2}\ar[ul]|-\circ&&
C_0(\Z^n) \ar[ul]|-\circ \ar[ll]|-\circ_{\Sigma \kappa_1}
}\nonumber
\end{align}
A sequence of distinguished triangles of this type will be called a \emph{tower}. The tower \eqref{zntower} in $KK^{\Z^n}$ is the higher rank analogue of the distinguished triangle \eqref{futt}. The tower \eqref{zntower} can be generalized to contain any coefficient ring.

To find a better description of the morphisms $\kappa_i$ let us recall the notion of a Koszul complex. Let $R$ denote a commutative ring and $E$ an $R$-module. For simplicity we will assume that $E$ is free and finitely generated, let us say of rank $N$. For an element $v\in\Hom_R(E,R)$, the Koszul complex of $E$ with respect to $v$ is the complex
\[0\to  \wedge^{N}E \xrightarrow{\partial_1}\wedge^{N-1}E\xrightarrow{\partial_2}\ldots \xrightarrow{\partial _{N-2}} \wedge^{2}E\xrightarrow{\partial_{N-1}} E\xrightarrow{v} R\to 0,\]
where each $\partial _k$ is defined as interior multiplication by $v$. Since we have assumed $E$ to be free, we may write $\nu=\sum \nu_i\e_i^*$ for some $\nu_1,\nu_2,\ldots \nu_N\in R$ and the dual basis $\e_i^*$ of a basis $\e_i$, $i=1,2,\ldots N$, of $E$. If the sequence $\nu_1,\nu_2,\ldots, \nu_N$ is a regular sequence the Koszul complex is exact except at $R$. The cohomology of the Koszul complex is in this case $R/v(E)$ at $R$. See more in \cite{eisen}.

The Koszul complex of interest to us is constructed from the module $E:=R(T)^n$ over the representation ring of the torus $T$ which has the following form:
\[R(T)\cong \Z[t_1^{\pm1}, \ldots, t_n^{\pm1}].\] 
Observe that Baaj-Skandalis duality and the Universal Coefficient Theorem implies that 
\begin{align*}
KK^{\Z^n}(\C^{k_i}\otimes C_0(\Z^n),\C^{k_{i+1}}\otimes C_0(\Z^n))&\cong KK^T(\C^{k_i},\C^{k_{i+1}})\\
&\cong \Hom_{R(T)}(R(T)^{k_i},R(T)^{k_{i+1}}).
\end{align*}
We have that $R(T)^{k_i}\cong \wedge^{n-i+1}E$ so the lower row in \eqref{zntower} have the right ranks for coinciding with a Koszul complex. Let $f_i\in \Hom_{R(T)}( \wedge^{n-i+1}E, \wedge^{n-i}E)$ denote the image of $\kappa_i$ under the isomorphisms above. To simplify notations, we will by $(e_i)_{i=1}^n$ denote the $R(T)$-basis of $E$ coming from the isomorphism $E\cong R(T)\otimes_\Z \Z^n$ and by $(e_i^*)_{i=1}^n$ denote the dual basis.

\begin{sats}
Under the isomorphisms $R(T)^{k_i}\cong \wedge^{n-i+1}E$ the mappings $f_i$ coincide with interior multiplication by the element $v:=\sum_{i=1}^n (1-t_i)\e_i^*$. Therefore the sequence 
\[0\to  \wedge^{n}E \xrightarrow{f_1}\wedge^{n-1}E\xrightarrow{f_2}\ldots \xrightarrow{f_{n-2}} \wedge^{2}E\xrightarrow{f_{n-1}} E\xrightarrow{f_n} R(T)\to 0\]
defines a complex isomorphic to the Koszul complex of $E$ whose cohomology at $R(T)$ is $\Z$.
\end{sats}

\begin{proof}
While both $f_i$ and the mapping defined by interior multiplication by $v$ are $R(T)$-linear it is sufficient to prove that $f_i(u)=v\neg u$ for elements of the form  $u=\e_{m_1}\wedge \cdots \wedge \e_{m_{n-i+1}}\in \wedge^{n-i+1}E$, where $m_1, \ldots , m_{n-i+1}\in\{1, 2, \ldots , n\}$. Let $(m_p)_{p=n-i+1}^n$ be an enumeration of all $j=1,2,\ldots ,n $ such that $j\notin (m_p)_{p=1}^{n-i+1}$. If we view $\Z^n$ as a subset of $\R^n$ we can define $X_u\subseteq \tilde{X}_i$ as the open face in $\R^n$ spanned by the vectors $\e_{m_{n-i+1}}, \e_{m_{n-i+2}}, \ldots \e_{m_n}$.

Under the isomorphism $\wedge^{n-i+1}E\cong K_{i-1}(\C^{k_i}\otimes \Sigma^{i-1}C_0(\Z^n))$ the element $u$ corresponds to a $K$-theory class on $\tilde{X}_i$ which is trivial except on the face $X_u$. Therefore there exists sequences of numbers $(a_j)_{j=1}^{n-i+1} , (b_j)_{j=1}^{n-i+1}\subseteq  \Z$ such that 
\[f_i(u)=\sum_{j=1}^{n-i+1} (a_j+b_jt_{j})\e_{m_j}\neg u.\]
If $j=1,2 \ldots, n-i+1$, we will let $X_{u,j}$ denote the open face spanned by the vectors $\e_{m_j}, \e_{m_{n-i+1}}, \e_{m_{n-i+2}}, \ldots \e_{m_n}$. It follows from restricting to $X_{u,j}$ that $a_j=1$ since Bott periodicity implies that the index mapping $K_{i-1}(C_0(X_{u}))\to K_i(C_0(X_{u,j}))$ is an isomorphism. In a similar fashion it follows that $b_i=-1$. 

While $v(E)$ is the ideal generated by the regular sequence $1-t_1,1-t_2, \ldots, 1-t_n$, the cohomology of the Koszul complex is $R(T)/v(E)=\Z$ and the quotient mapping $R(T)\to \Z$ coincides with the augmentation mapping. 
\end{proof}

Consider the tower Baaj-Skandalis dual to \eqref{zntower}. Given $A,B\in KK^{T}$ we can apply the homological functor $KK^T(A,-\otimes_{min} B)$ to this tower. This functor is only homological on the bootstrap category if $B$ is not exact, but all objects in the tower Baaj-Skandalis dual to \eqref{zntower} are in the bootstrap category. The lowest row of the corresponding tower in the category of $R(T)$-modules is a Koszul complex:
\begin{align}
\label{kcalc}
0\to  \wedge^{n}\Z^n\otimes  KK_*^T(A,B) &\xrightarrow{v_A\neg}\wedge^{n-1}\Z^n\otimes  KK_*^T(A,B)\xrightarrow{v_A\neg}\ldots\\
\ldots&\xrightarrow{v_A\neg} \Z^n\otimes  KK_*^T(A,B)\xrightarrow{v_A\neg} KK_*^T(A,B) \to 0\nonumber
\end{align}
where 
\[v_A:=\sum_{i=1}^n (1-\beta_{i*})\e_i^*\in \Hom_{R(T)}(KK^T_*(A,B)^n,KK^T_*(A,B))\] 
and $(\beta_i)_{i=1}^n$ are the commuting equivariant automorphisms of $A$ that are Baaj-Skandalis to the $\Z^n$-action on $B\rtimes_rT$. The cohomology of this Koszul complex can be calculated from $KK_*^T(A,B)$. We will return to this subject in the next section in the more general case of Hodgkin-Lie groups and explain this procedure further.

\section{The generalized Pimsner-Voiculescu-towers}

As mentioned in the introduction, the representation ring $R(T)$ is free over $R(G)$ when $G$ is a Hodgkin-Lie group, so the step to coactions of a compact Hodgkin-Lie group will not be too large. We will throughout this section assume that $G$ is a compact Hodgkin-Lie group of rank $n$ with maximal torus $T$. Recall that a group satisfies the Hodgkin condition if it is connected and the fundamental group is torsion-free. 

The embedding $T\subseteq G$ induces a restriction functor $KK^{\hat{T}}\to KK^{\hat{G}}$. Using the isomorphism $\hat{T}\cong \Z^n$, the tower \eqref{zntower} can be restricted to a $KK^{\hat{G}}$-tower:
\begin{align*}
\xymatrix@C=0.2em@R=2.71em{
& \Sigma^nC^*(T) \ar[rr]&&
C_0(Y_{n-1}) \ar[rr] \ar[dl] &&
\cdots \ar[dl]\\
&& \C^{n}\otimes \Sigma^{n-1}C^*(T) \ar[ul]|-\circ_{\Sigma^n\kappa_n}&&
\C^{k_{n-1}}\otimes \Sigma^{n-2} C^*(T) \ar[ul]|-\circ \ar[ll]|-\circ_{\Sigma^{n-1}\kappa_{n-1}}&&
\quad \cdots\quad \ar[ul]|-\circ \ar[ll]|-\circ}\\
\xymatrix@C=0.2em@R=2.71em{
\cdots \ar[rr] &&
C_0(Y_2) \ar[rr] \ar[dl] &&
C_0(Y_1) \ar[dl]\ar[rr]&&
C_0(\R^n) \ar[dl] \\
& \cdots \quad&& 
\C^{n}\otimes \Sigma C^*(T)\ar[ll]|-\circ_{\Sigma^2\kappa_2}\ar[ul]|-\circ&&
C^*(T) \ar[ul]|-\circ \ar[ll]|-\circ_{\Sigma \kappa_1}
}
\end{align*}
In order to work with this $KK^{\hat{G}}$-tower we need to describe the terms $C^*(T)$ in the second row.

\begin{lem}
\label{ishom}
If $G$ is a compact Hodgkin-Lie group with Weyl group of order $w$ there is an isomorphism 
\[C^*(T)\cong \C^{w}\otimes C^*(G)\quad \mbox{in} \quad KK^{\hat{G}}.\]
\end{lem}

Observe that the condition on $G$ to be a Hodgkin group is equivalent to $\hat{G}$ being a torsion-free quantum group in the sense of Meyer, see \cite{mehom}. The torsion-free quantum groups are the only non-classical discrete quantum groups for which there is a general formulation of the Baum-Connes property in terms of triangulated categories. In \cite{nemetva}, coactions of compact non-Hodgkin Lie groups were considered and the "torsion" turned out to be the torsion elements of $H^2(G,S^1)$. The less precise statement $C(G/T)\cong \C^k$ in $KK^G$ for some $k$ is stated and proved in \cite{nemetva}. An explicit calculation that $k=|W|$ can be found in \cite{stei}. We will review the conceptually important  part of the proof of a Proposition in \cite{nemetva} which proves Lemma \ref{ishom} aside from the calculation of $k$.

\begin{proof} 
By \cite{stei}, the representation ring $R(T)$ is free of rank $w$ over the representation ring $R(G)$ if $\pi_1(G)$ is torsion-free.  If we let $\mathcal{S}$ denote the localizing subcategory of $KK^G$ generated by $\C$ and $C(G/T)$, Lemma $11$ of \cite{neme} states that for $A\in \mathcal{S}$ the natural homomorphism 
\[R(T)\otimes _{R(G)}KK^G(A,\C)\to KK^T(A,\C)\]
is an isomorphism. Thus the representable functor on $\mathcal{S}$
\[A\to KK^G(A,\C^{w})\cong R(T)\otimes _{R(G)}KK^G(A,\C)\] 
coincides with the representable functor 
\[A\to KK^G(A,C(G/T))\cong KK^T(A,\C).\] 
The last isomorphism exists as a consequence of the fact that the induction functor $Ind^G_T$ is the right adjoint of the restriction functor from $G$ to $T$. So the Yoneda lemma implies that $C(G/T)\cong \C^{w}$ in $\mathcal{S}$ and therefore in $KK^G$. Applying Baaj-Skandalis duality it follows that there is an equivariant $KK$-isomorphism $C^*(T)\cong \C^{w}\otimes C^*(G)$.
\end{proof}

Using Lemma \ref{ishom} the tower \eqref{zntower} takes the form:
\begin{align}
\label{fpvtower}
\xymatrix@C=0.2em@R=2.71em{
& \Sigma^n\C^w\otimes C^*(G) \ar[rr]&&
C_0(Y_{n-1}) \ar[rr] \ar[dl] &&
\cdots \\
&& \Sigma^{n-1}\C^{wn}\otimes C^*(G) \ar[ul]|-\circ&&
\;\cdots \ar[ul]|-\circ \ar[ll]|-\circ
}\\
\xymatrix@C=0.2em@R=2.71em{
\cdots \ar[rr] &&
C_0(Y_2) \ar[rr] \ar[dl] &&
C_0(Y_1) \ar[dl]\ar[rr]&&
\Sigma^n\C \ar[dl] \\
& \cdots \quad&&  \Sigma \C^{wn}\otimes C^*(G)\ar[ll]|-\circ\ar[ul]|-\circ&&
\C^w\otimes C^*(G) \ar[ul]|-\circ
 \ar[ll]|-\circ
}\nonumber
\end{align}
We will call this $KK^{\hat{G}}$-tower the fundamental $G-$PV-tower. The dual fundamental $G-$PV-tower is defined to be the $KK^{G}$-tower which is Baaj-Skandalis dual to the fundamental $G-$PV-tower:
\begin{align}
\label{dfpvtower}
\xymatrix@C=0.2em@R=2.71em{
& \Sigma^n\C^w \ar[rr]&&
D_{n-1} \ar[rr] \ar[dl] &&
\cdots \ar[dl]\\
&& \Sigma^{n-1}\C^{wn} \ar[ul]|-\circ&&
\Sigma^{n-2} \C^{wk_{n-1}} \ar[ul]|-\circ
 \ar[ll]|-\circ&&
\cdots \ar[ul]|-\circ \ar[ll]|-\circ
}\\
\xymatrix@C=0.2em@R=2.71em{
\cdots \ar[rr] &&
D_2 \ar[rr] \ar[dl] &&
D_1 \ar[dl]\ar[rr]&&
\Sigma^nC(G) \ar[dl] \\
& \cdots \quad&&  \Sigma \C^{wn}\ar[ll]|-\circ\ar[ul]|-\circ&&
\C^w \ar[ul]|-\circ
 \ar[ll]|-\circ
}\nonumber
\end{align}
where $D_i:=C_0(Y_i)\rtimes_r\hat{G}$. 

As mentioned above, if $A$ is a $G-C^*$-algebra, the trivial coaction of $G$ on $A$ makes $A$ into a Yetter-Drinfeld algebra. This follows from that $C(G)$ is commutative so we can extend a $G$-action via the $D(G)$-equivariant $*$-monomorphism $C(G)\to \mathcal{M}(C_0(D(G)))$. Clearly, a $G$-equivariant mapping is equivariant in this new $D(G)$-action. Furthermore, since mapping cones does not depend on the action, the trivial extension of a $G$-action to a $D(G)$-action is functorial and respects mapping cones. The following proposition follows from universality.

\begin{prop}
\label{doubleext}
If $G$ is a locally compact group, the functor mapping a $G-C^*$-algebra to a $G$-Yetter-Drinfeld algebra with trivial $\hat{G}$-action defines a triangulated functor $KK^G\to KK^{D(G)}$.
\end{prop}

Using the triangulated functor of Proposition \ref{doubleext}, we may consider the tower \eqref{dfpvtower} as a tower in $KK^{D(G)}$. Applying a crossed product by $G$ we obtain that also the tower \eqref{fpvtower} is a tower in $KK^{D(G)}$. For a $C^*$-algebra $B$ we will use the notation $t(B)$ for the $\hat{G}-C^*$-algebra with trivial coaction, or in the context of $G-C^*$-algebras $t(B)$ will denote the $G-C^*$-algebra with trivial action. Let us state and prove the corresponding version of \eqref{crobra} in a simple case of a braided tensor product over $G$ with $C(G)$, a more general proof can be found in \cite{voigt}.

\begin{lem}
\label{tenscross}
When $B$ has a continuous $G$-action, there is a $\hat{G}$-equivariant Morita equivalence 
\[(C(G)\otimes B)\rtimes_r G\sim_M t(B).\]
\end{lem}

\begin{proof}
By Baaj-Skandalis duality, it suffices to prove that there is a $\hat{G}$-equivariant isomorphism $(C(G)\otimes B)\rtimes_r G\cong (C(G)\rtimes_r G)\otimes t(B)$. Denote the $G$-action on $B$ by $\beta$ and define the equivariant mapping $\phi_0:L^1(G,C(G,B))\to (C(G)\rtimes_r G)\otimes t(B)$ by setting 
\[\phi_0(f)(g_1,g_2):=\beta_{g_1^{-1}}f(g_1,g_2).\]
The linear mapping $\phi_0$ is a $*$-homomorphism when $L^1(G,C(G,B))$ is equipped with the convolution twisted by the $G$-action on $C(G)\otimes B$. It is straightforward to verify that $\phi_0$ is bounded in $C^*$-norm so we can define $\phi:(C(G)\otimes B)\rtimes_rG\to (C(G)\rtimes_r G)\otimes B$ by continuity. The $*$-homomorphism $\phi$ is an equivariant isomorphism since an inverse can be constructed by extending 
\[\phi^{-1}(f\otimes b)(g_1,g_2):=f(g_1,g_2)\beta_{g_1}(b)\] 
to a $*$-homomorphism $\phi^{-1}:(C(G)\rtimes_r G)\otimes t(B)\to (C(G)\otimes B)\rtimes_rG$.
\end{proof}

\begin{sats}[The Pimsner-Voiculescu tower]
\label{gpv}
Let $G$ be a compact Hodgkin-Lie group of rank $n$ and Weyl group of order $w$. For any separable $\hat{G}-C^*$-algebra $A$ there is a $KK^{\hat{G}}$-tower 
\begin{align}
\label{pwtower}
\xymatrix@C=0.2em@R=2.71em{
& \C^w\otimes A\ar[rr]&&
\Sigma^nD_{n-1}(A) \ar[rr] \ar[dl] &&
\Sigma^nD_{n-2}(A) \ar[rr] \ar[dl] &&
\cdots \\
&& \Sigma\C^{wn}\otimes A \ar[ul]|-\circ&&
\Sigma^{2}\C^{wk_{n-1}}\otimes A \ar[ul]|-\circ\ar[ll]|-\circ&&
\cdots\ar[ul]|-\circ\ar[ll]|-\circ
}\\
\xymatrix@C=0.2em@R=2.71em{
\cdots \ar[rr] &&
\Sigma^nD_2(A) \ar[dl]\ar[rr]&&
\Sigma^nD_1( A) \ar[dl]\ar[rr]&&
t(A\rtimes_r\hat{G})\ar[dl] \\
& \cdots &&
\Sigma^{n-1}\C^{wk_2}\otimes A \ar[ul]|-\circ \ar[ll]|-\circ &&
\Sigma^n\C^w\otimes A \ar[ul]|-\circ \ar[ll]|-\circ
}\nonumber
\end{align}
where $D_i(A):=(C_0(Y_i)\otimes \Ko(L^2(G)))\boxtimes_G (A\rtimes_r \hat{G})$ and is equipped with the $\hat{G}$-action induced from the diagonal $\hat{G}$-action on $C_0(Y_i)\otimes \Ko(L^2(G))$.
\end{sats}

Observe that the $D(G)$-actions on the $C^*$-algebras $C_0(Y_i)\otimes \Ko(L^2(G))$ is defined to come from those on their Baaj-Skandalis duals $C_0(Y_i)\rtimes_r \hat{G}$, which are $D(G)-C^*$-algebras in the dual $G$-actions on the crossed products and the trivial $\hat{G}$-actions. So in general, $D_i(A)$ is not the tensor product of $C_0(Y_i)\otimes \Ko(L^2(G))$ and $A\rtimes _r\hat{G}$.

\begin{proof}
By Lemma \ref{tenscross} the $\hat{G}-C^*$-algebra $A$ admits the equivariant Morita equivalence:
\begin{equation}
\label{firmor}
(C(G)\otimes (A\rtimes_r \hat{G}))\rtimes_rG\sim_M t(A\rtimes_r \hat{G}).
\end{equation}
Furthermore, the isomorphism of equation \eqref{crobra} holds for braided tensor products over $G$ so while the $\hat{G}$-actions on $D_i=C_0(Y_i)\rtimes_r\hat{G}$ are trivial there are equivariant isomorphisms
\begin{align}
\label{seciso}
(D_i\otimes (A\rtimes_r \hat{G}))\rtimes_rG&\cong((C_0(Y_i)\rtimes_r \hat{G})\boxtimes_G (A\rtimes_r \hat{G}))\rtimes_rG\cong\\
&\cong(C_0(Y_i)\otimes \Ko(L^2(G)))\boxtimes_G (A\rtimes_r \hat{G}).\nonumber
\end{align}
Thus if we tensor the dual fundamental $G-$PV-tower \eqref{dfpvtower} by the $G-C^*$-algebra $A\rtimes_r \hat{G}$ we obtain a new $KK^{G}$-tower which becomes the Pimsner-Voiculescu tower of $A$ after applying Baaj-Skandalis duality, using the Morita equivalence \eqref{firmor} and the isomorphisms \eqref{seciso}.
\end{proof}

The Pimsner-Voiculescu tower \eqref{pwtower} is the generalization of the resolution in \eqref{kcalc} to compact Hodgkin-Lie groups. Applying the cohomological functor $KK(-,B)$ to the Pimsner-Voiculescu tower we obtain a similar resolution of $KK_*(A\rtimes_r\hat{G},B)$ in terms of $KK_*(A,B)$ as in \eqref{kcalc}. Similarly, the homological functor $KK(B,-)$ applied to the Pimsner-Voiculescu tower gives a resolution of $KK(B,A\rtimes \hat{G})$ in terms of $KK(B,A)$. Observe that since $A$ has a $\hat{G}$-action, the groups $KK(\C^w\otimes A,B)$ and $KK(B,\C^w\otimes A)$ will always have an $R(G)$-module structure and since $R(T)$ is free over $R(G)$ also an $R(T)$-module structure.

As an example of this, we will use the Pimsner-Voiculescu tower to calculate the $K$-theory of the homogeneous space $G/H$ when $H\subseteq G$ is a Lie subgroup. More generally, this technique can be used to calculate $K_*(A\rtimes_r \hat{G})$ for any $\hat{G}-C^*$-algebra $A$ when one knows $K_*(A)$ and its $R(G)$-module structure coming from the Julg isomorphism $K_*(A)\cong K_*^G(A\rtimes_r\hat{G})$. To calculate $K^*(G/H)$, consider the $C^*$-algebra $A:=C^*(H)$ equipped with the $\hat{G}$-action induced from the natural mapping $C^*(H)\to \mathcal{M}(C^*(G))$. Green's imprimitivity theorem implies that $C^*(H)\rtimes \hat{G}$ is $KK$-equivalent with $C(G/H)$. Thus, if we take the $K$-theory of the Pimsner-Voiculescu tower of $C^*(H)$ we obtain a tower of abelian groups of the form 
\begin{align}
\label{gmodh}
\xymatrix@C=0.2em@R=2.71em{
& R(T)\otimes_{R(G)} R(H)\ar[rr]&&
K_{*-n}(D_{n-1}(C^*(H))) \ar[rr] \ar[dl] &&
\cdots \\
&& \Sigma R(T)^n\otimes_{R(G)} R(H) \ar[ul]|-\circ _{\Sigma v\otimes1}&&
\cdots\ar[ul]|-\circ\ar[ll]|-\circ_{\Sigma v\otimes 1}
}\\
\xymatrix@C=0.2em@R=2.71em{
\cdots \ar[rr] &&
K_{*-n}(D_{1}(C^*(H)))\ar[dl]\ar[rr]&&
K^*(G/H)\ar[dl] \\
& \cdots &&
\Sigma^n R(T)\otimes_{R(G)} R(H) \ar[ul]|-\circ \ar[ll]|-\circ _{\Sigma v\otimes 1}
}\nonumber
\end{align}
We use $\Sigma$ to denote degree shift in the category of $\Z/2\Z$-graded abelian groups. Here we have used that $R(T)$ is a free $R(G)$-module of rank $w$ so $K_*(\C^w\otimes C^*(H))\cong R(T)\otimes_{R(G)} R(H)$. Thus the lowest row is the tensor product of $R(H)$ with the Koszul complex of $R(T)$ that is associated with the regular sequence $1-t_1,1-t_2, \ldots ,1-t_n$ under the isomorphism $R(T)\cong \Z[t_1^{\pm1},t_2^{\pm 1},\ldots, t_n^{\pm 1}]$. 

If we restrict our attention to simple compact Lie groups we can perform an explicit calculation of all the groups in \eqref{gmodh}. Assume that $G=G_n$ is a simple compact Hodgkin-Lie group in the classical $A,B,C$- or $D$-series of rank $n$ and assume that $H=G_k\subseteq G_n$ is a simple simply connected compact Lie group in the same classical serie being of rank $k<n$. We may take a maximal torus $T_n\subseteq G_n$ such that $T_k:=T_n\cap G_k$ is a maximal torus in $G_k$. In this case we may consider $R(T_k)$ as an ideal in $R(T_n)$ and $R(T_n)\otimes _{R(G_n)}R(G_k)\cong R(T_k)$ as $R(T_n)$-modules. Under the isomorphisms $R(T_k)\cong \Z[t_1^{\pm1},t_2^{\pm 1},\ldots, t_k^{\pm 1}]$ and $R(T_n)\cong \Z[t_1^{\pm1},t_2^{\pm 1},\ldots, t_n^{\pm 1}]$, the Koszul vector $v$ is identified with $\sum_{i=1}^k (1-t_i)\e_i^*\in \Hom(R(T_k)^n,R(T_k))$. Thus we arrive at the tower
\begin{align*}
\xymatrix@C=0.01em@R=2.71em{
& R(T_k)\ar[rr]&&
K_{*-n}(D_{n-1}(C^*(G_k)))) \ar[rr] \ar[dl] &&
\cdots \\
&& \Sigma \Z^n\otimes_\Z R(T_k) \ar[ul]|-\circ _{\Sigma\partial_n}&&
\wedge^2\Z^n\otimes_\Z R(T_k) \ar[ul]|-\circ\ar[ll]|-\circ_{\Sigma\partial_{n-1}}&&
\cdots\ar[ll]|-\circ
}\\
\xymatrix@C=0.2em@R=2.71em{
\cdots \ar[rr] &&
K_{*-n}(D_{1}(C^*(G_k)))\ar[dl]\ar[rr]&&
K^*(G_n/G_k)\ar[dl] \\
& \cdots &&
\Sigma^n \wedge^n\Z^n\otimes_\Z R(T_k) \ar[ul]|-\circ \ar[ll]|-\circ _{\Sigma \partial_1}
}
\end{align*}
Let us use the notation $E^*$ for the complex $\wedge^{n-*}\Z^n\otimes R(T_k)$ equipped with the Koszul differential from the vector $\sum_{i=1}^k (1-t_i)\e_i^*$ which we as above denote by $\partial_l:E^{l-1}\to E^l$. After some simpler calculations in this Koszul complex we arrive at the conclusion that 
\[K_{*-n}(D_{l}(C^*(G_k)))\cong \ker(\partial_{l+1})\oplus \bigoplus_{j=l+2}^{n+1}\Sigma^{n-j}H^j(E^*).\]
Hence we obtain the isomorphism $K^*(G_n/G_k)\cong \bigoplus_{j=1}^{n+1}\Sigma^{n-j}H^j(E^*)$. These cohomology groups are calculated in Corollary $17.10$ of \cite{eisen} and $H^j(E^*)$ is a free group of rank $k(j):=(n-k)!/(n-j)!(n-j-k)!$ if $0\leq j\leq n-k$ and $0$ otherwise. Therefore 
\[K^{*}(G_n/G_k)\cong \bigoplus_{j=0}^{n-k}\Sigma^{n-j}\Z^{k(j)}=\Z^{2^{n-k-1}}\oplus \Sigma \Z^{2^{n-k-1}}.\]

\begin{sats}[The dual Pimsner-Voiculescu tower]
\label{dgpv}
Under the assumptions of Theorem \ref{gpv} there is a $KK^{G}$-tower
\begin{align}
\label{dualpwtower}
\xymatrix@C=0.2em@R=2.71em{
& \C^w\otimes t(A)\ar[rr]&&
\Sigma^n\tilde{D}_{n-1}(A) \ar[rr] \ar[dl] &&
\;\cdots \\
&& \Sigma\C^{wn}\otimes t(A) \ar[ul]|-\circ&&
\Sigma^{2}\C^{wk_{n-1}}\otimes t(A) \ar[ul]|-\circ\ar[ll]|-\circ&&
\;\;\cdots\ar[ll]|-\circ
}\\
\xymatrix@C=0.2em@R=2.71em{
\cdots \ar[rr] &&
\Sigma^n\tilde{D}_2(A) \ar[dl]\ar[rr]&&
\Sigma^n\tilde{D}_1(A) \ar[dl]\ar[rr]&&
A\rtimes_r\hat{G}\ar[dl] \\
&\cdots &&
\Sigma^{n-1}\C^{wk_2}\otimes t(A) \ar[ul]|-\circ \ar[ll]|-\circ &&
\Sigma^n\C^w\otimes t(A) \ar[ul]|-\circ \ar[ll]|-\circ
}\nonumber
\end{align}
where $\tilde{D}_i(A):=D_i\boxtimes_{\hat{G}} A$.
\end{sats}

For a homological functor $F:KK^{\hat{G}}\to Ab$, the dual Pimsner-Voiculescu tower of $A$ allows us to calculate $F(A)$ in terms of the objects $F(C^*_r(G)\otimes t(A))$. As we shall see below, $\hat{G}-C^*$-algebras of the form $C^*_r(G)\otimes t(A)$ behaves similarly to \emph{proper} actions. Compare this result to Theorem $4.4$ of \cite{mehomett}.

\begin{proof}
Consider the braided tensor product by $\Sigma^n A$ and the tower \eqref{fpvtower}:
\begin{align*}
\xymatrix@C=0.2em@R=2.71em{
&\C^w\otimes C^*(G)\boxtimes_{\hat{G}}  A \ar[rr]&&
\Sigma^n C_0(Y_{n-1})\boxtimes_{\hat{G}}  A \ar[rr] \ar[dl] &&
\cdots \\
&& \Sigma\C^{wn}\otimes C^*(G)\boxtimes_{\hat{G}}  A \ar[ul]|-\circ&&
\;\cdots \ar[ul]|-\circ \ar[ll]|-\circ
}\\
\xymatrix@C=0.2em@R=2.71em{
\cdots \ar[rr] &&
\Sigma^n C_0(Y_1)\boxtimes_{\hat{G}}  A \ar[dl]\ar[rr]&&
A \ar[dl] \\
& \cdots \quad&&  
\Sigma^n\C^w\otimes C^*(G)\boxtimes_{\hat{G}}  A \ar[ul]|-\circ
 \ar[ll]|-\circ
}
\end{align*}
Taking crossed product between this tower and $\hat{G}$ implies the Theorem since the following equivariant Morita equivalences follows from \eqref{crobra} 
\begin{align*}
(C^*(G)\boxtimes_{\hat{G}}  A)\rtimes_r\hat{G}&\sim_M t(A) \quad\mbox{and}\\
(C_0(Y_i)\boxtimes_{\hat{G}}  A)\rtimes_r\hat{G}&\sim_M (C_0(Y_i)\rtimes_r\hat{G})\boxtimes_{\hat{G}} A=D_i\boxtimes_{\hat{G}} A.
\end{align*}
\end{proof}

One of the main motivations behind this paper was to give a precise description of the Baum-Connes property of duals of Hodgkin-Lie groups. The Baum-Connes property for coactions of compact Lie groups was given meaning to and was proved to hold in \cite{nemetva}. More generally, this fits into the program of generalizing the Baum-Connes property to quantum groups. So far, it is not known what a suitable property the Baum-Connes property should be for a general locally compact quantum group. For discrete quantum groups which are torsion-free, in the sense of \cite{mehom}, there is a formulation and as mentioned above duals of compact Hodgkin-Lie groups are torsion-free. 

The problem that arises when one tries to define the Baum-Connes assembly mapping for a quantum group is that there is no natural notion of a proper action and there are in general too many quantum homogeneous spaces. It is much easier to generalize certain notions of free actions than proper actions of a quantum group by just saying that an action of a discrete quantum group $\Gamma$ on a $C^*$-algebra $A$ is \emph{truly} free if there is a $C^*$-algebra $A_0$ and an equivariant $*$-isomorphism $A\cong A_0\otimes_{min} C_0(\Gamma)$ with $\Gamma$ only acting on the second leg. In the case of a group, there are many free actions that are not truly free but this stronger notion of a free action will suffice for our purposes. 

Restricting one's attention to generalizing the Baum-Connes property of the simpler class of torsion-free discrete groups to the quantum setting, when proper actions are free, Meyer introduced a class of quantum groups known as torsion-free in \cite{mehom}. Following \cite{mehom}, we say that a discrete quantum group $\Gamma$ is torsion-free if every coaction of the compact quantum group $\hat{\Gamma}$ on a finite-dimensional $C^*$-algebra is Morita equivalent to a trivial coaction on a direct sum of $\C$:s. This fact implies that any finite-dimensional projective representation of the dual compact quantum group is equivalent to a representation. If $\Gamma$ is a discrete group, coactions of the dual compact quantum group on finite-dimensional $C^*$-algebras that are not Morita equivalent to a trivial coaction on a direct sum of $\C$:s correspond to finite subgroups so a discrete group is torsion-free if and only if it is torsion-free in the sense of \cite{mehom}.

For a torsion-free quantum group a proper action should correspond to a free action. Under Baaj-Skandalis duality, a truly free $\Gamma-C^*$-algebra corresponds to a trivial $\hat{\Gamma}$-action. Let $\mathcal{CI}_{\hat{\Gamma}}$ denote the image of $t:KK\to KK^{\hat{\Gamma}}$. The triangulated category $\langle \mathcal{CI}_{\hat{\Gamma}}\rangle $ is defined as the localizing subcategory generated by $\mathcal{CI}_{\hat{\Gamma}}$. Following the formulation of \cite{mehom}, $\Gamma$ is said to satisfy the strong Baum-Connes property if the embedding of triangulated categories $\langle\mathcal{CI}_{\hat{\Gamma}}\rangle\to KK^{\hat{\Gamma}}$ is essentially surjective. The strong Baum-Connes property of $\Gamma$ is equivalent to that any $\Gamma-C^*$-algebra is in the localizing category generated by all truly free actions. So regardless of what notion of a proper action we choose, the strong Baum-Connes conjecture will imply that the localizing category generated by all such proper actions will be $KK^\Gamma$. The quantum group is said to satisfy the Baum-Connes property if the same statement holds after localizing with respect to the kernel of equivariant $K$-theory. 

In \cite{nemetva} the Baum-Connes property was formulated in the slightly more general setting of duals of compact Lie groups. The finite-dimensional projective representations of a compact Lie group $G$ correspond to the torsion classes of $H^2(G,S^1)$, which can be thought of as the torsion of $\hat{G}$. When $G$ is Hodgkin, $H^2(G,S^1)$ is torsion-free so $\hat{G}$ is torsion-free. In this case a "proper" action is an object of the additive category generated by $\hat{G}$-algebras that are Baaj-Skandalis dual to $A_0\otimes \C_\omega$, with $\C_\omega$ denoting the endomorphisms of a projective representation $\omega$ and $A_0$ having trivial $G$-action. So the substitute in the setting of \cite{nemetva} for proper actions is the category of tensor products between Baaj-Skandalis duals of coactions on finite-dimensional $C^*$-algebras and trivial actions, just as the truly free actions form a substitute for proper actions of torsion-free quantum groups. The Baum-Connes property of coactions of a compact Hodgkin-Lie group is a direct consequence of Theorem \ref{dgpv}. The method of proof of Proposition $2.1$ of \cite{nemetva} can be used to generalize both Theorem \ref{gpv} and Theorem \ref{dgpv} to arbitrary compact Lie group. 

Finally, let us mention a promising generalization of Theorem \ref{dgpv} to Woronowicz deformations. It was proved in \cite{voigt} that the compact quantum group $SU_q(2)$ satisfies that $C(SU_q(2)/T)$ is $KK^{D(SU_q(2))}$-isomorphic to $\C^2$ for $q\in ]0,1[$. So if we apply the induction functor $Ind^{SU_q(2)}_T:KK^T\to KK^{SU_q(2)}$ to the distinguished triangle Baaj-Skandalis dual to \eqref{futt} and use the isomorphism of Nest-Voigt we arrive at the distinguished triangle in $KK^{D(SU_q(2))}$:
\[
\xymatrix@C=0.2em@R=2.71em{
\C^2\ar[rr] & & \C^2 \ar[dl]\\
  &C(SU_q(2)). \ar[ul]|-\circ
}
\]
Using the technique from the proof of Theorem \ref{dgpv} any $A\in KK^{\widehat{SU_q(2)}}$ fits into a distinguished triangle 
\[
\xymatrix@C=0.2em@R=2.71em{
\C^2\otimes t(A)\ar[rr] & & \C^2\otimes t(A) \ar[dl]\\
  &A\rtimes\widehat{SU_q(2)}. \ar[ul]|-\circ
}
\]
This distinguished triangle gives an alternative proof of the strong Baum-Connes property for $\widehat{SU_q(2)}$, a result first proved in \cite{voigtfree}. The interesting part about this proof is that it only relies on the isomorphism $C(G_q/T)\cong \C^w$ in $KK^{D(G_q)}$. So if such an isomorphism exists for a simply connected semi-simple compact Lie group $G$, the strong Baum-Connes conjecture holds for $\hat{G}_q$, the quantum dual of the Woronowicz deformation of $G$. To formulate the Baum-Connes property for $\hat{G}_q$ we must of course know that it is torsion-free, a statement proved in \cite{voigtfree} for $G=SU(2)$ and the general case was proved in \cite{gofftwist}. Another striking application of such an isomorphism is that the method above for calculating $K$-theory of homogeneous spaces can be generalized to classical quantum homogeneous spaces of the Woronowicz deformations.

\end{document}